\documentclass[a4paper,10pt]{article}
\usepackage{amssymb,amsmath,amsthm,latexsym}
\usepackage{color}
\usepackage{textcomp}
\usepackage[active]{srcltx}
\usepackage{colortbl}

\usepackage{color}
\usepackage{mathrsfs}
\usepackage{graphicx}
\usepackage{lineno}
\usepackage{adjustbox}
\usepackage{rotating}

\usepackage[multiple]{footmisc}

\usepackage{fancyhdr}

\newcommand{\keywords}[1]{%
  \small
  \textbf{\textit{Keywords---}} #1
}

\newcommand{\mscclass}[1]{\par\noindent\textbf{MSC:} #1}

\theoremstyle{definition}

\theoremstyle{remark}

\numberwithin{equation}{section}

\title{Small degree Salem numbers with trace $-3$}

\author{J.-M. Sac-\'Ep\'ee\footnote{J.-M. Sac-\'Ep\'ee, IECL, Université de Lorraine,
 France, jean-marc.sac-epee@univ-lorraine.fr}}

\begin{document}

\maketitle
\newcommand{\D}{\mathbb{D}}
\newcommand{\C}{\mathbb{C}}
\newcommand{\N}{\mathbb{N}}
\newcommand{\R}{\mathbb{R}}
\newcommand{\Z}{\mathbb{Z}}
\newcommand{\dist}{\operatorname{dist}}

\renewcommand{\qedsymbol}{$\blacksquare$}

\makeatletter
\def\blfootnote{\xdef\@thefnmark{}\@footnotetext}
\makeatother

\begin{abstract}
The contribution of this work is to provide tables of Salem numbers with trace $-3$ and small degrees, namely degrees $2d = 34, 36, 38$, and $40$. The implemented method also generates a list of totally positive polynomials of degrees $d = 17, 18, 19$, and $20$, with trace $2d-3$,  a single root strictly greater than $4$, and all other roots strictly between $0$ and $4$.
\end{abstract}
\mscclass{11K16}


\section{Introduction} 

A Salem number is a real algebraic integer strictly greater than 1, whose conjugates all have a modulus less than or equal to 1, and at least one conjugate has a modulus equal to 1. The minimal polynomial P of a Salem number (referred to as a {\it Salem polynomial}) is self-reciprocal, i.e. $x^d P(1/x) = P (x)$, where $d$ is the  degree of $P$. A direct consequence of this is that a Salem number is conjugate to its inverse, and the rest of its conjugates are on the unit circle.

By definition, the degree of a Salem number is simply that of its minimal polynomial. This degree must be even and at least equal to 4. Similarly, the trace of a Salem number is that of its minimal polynomial, i.e. the sum of its roots. With $P(x) = \displaystyle\sum_{k=0}^d a_k x^k$ ($a_d = 1$), the trace of $P$ is equal to $-a_{d-1}$.

A very comprehensive survey devoted to Salem numbers is \cite{Smyth2015}, where the essential results and properties proved on the subject over the last few decades are presented and explained in detail.

As stated in \cite{McKeeSmyth2005}, examples of Salem numbers of any strictly positive trace are available since the largest real root of the polynomial $x^4 - nx^3 -(2n+1)x^2 -nx+1$  is a Salem number (of trace $n$) for all $ n \geq 1$.  It is also well known that  the polynomial $x^6-x^4-2x^3-x^2+1$  is the minimal polynomial of a Salem number of zero trace.

Things get much more complicated when it comes to finding Salem numbers with negative traces, and a great deal of work has been devoted to this.

In 2004, J. McKee and C.J. Smyth  provided Salem numbers with trace $-2$ \cite{McKeeSmyth2004}. The same two authors proved in 2005  that there are Salem numbers of all negative traces \cite{McKeeSmyth2005}. In 2011, J. McKee \cite{McKee2011} provided an explicit example of a Salem number of trace $-3$ (its minimal polynomial was of degree 54). Four years later, S. El Otmani, G. Rhin and J.-M. Sac-Épée  provided a Salem number with degree $34$ and  trace $-3$ \cite{ElOtmaniRhinSac-Epee2015}, expressing the conviction, motivated by the results of \cite{ElOtmaniMaulRhinSac-Epee2014}, that there is presumably no Salem number with trace $-3$ and degree $32$. In 2024, G. Cherubini and P. Yatsyna proved \cite{CherubiniYatnyna2024} that there exist Salem numbers with trace $-3$ and every even degree $\geq 34$, notably providing a table containing an example of a Salem number of trace $-3$ for each of the degrees from $34$ to $52$, as well as an example at degree $58$.

In the present work, our aim is to provide tables giving numerous examples of Salem numbers of trace $-3$ and small degrees, i.e. degrees between $34$ and $40$. The method on which our approach is based, namely integer linear optimization (since the variables of the problem are the \textit{integer} coefficients of the polynomials to be found), has already enabled us to provide a list of totally positive polynomials of degree $d =16$ and trace $2d-3 = 29$ \cite{ElOtmaniMaulRhinSac-Epee2014}, as well as the Salem number of trace $-3$ and degree $34$ previously mentioned  \cite{ElOtmaniRhinSac-Epee2015}. This method is explained in detail in the two articles cited above, so we will just summarize the main points for the reader's convenience.

\section{Outline of the method used to find Salem numbers}
\subsection{Relationship between Salem numbers and totally positive polynomials}
Consider a Salem polynomial $P$ with degree $2d$. By the change of variable $z = x + \displaystyle\frac{1}{x} + 2$, we obtain a monic integer polynomial $Q$ of degree $d$. This polynomial $Q$ is characterized by all its roots being real and positive: $d - 1$ roots are in $(0, 4)$, and only one root exceeds $4$. We then set ourselves the goal of finding totally positive monic polynomials of degree $d = 17, 18, 19, 20$ and trace $2d - 3 = 31, 33, 35, 37$, with $d-1$ roots in $(0,4)$ and only one root greater than $4$. In addition to its own interest, each of the polynomials $Q$ of degree $d$ found will immediately provide us with a Salem number of degree $2d$ via an inverse change of variable.

\subsection{Algorithmic approach}

The initial phase consists of locating the zeros of such a polynomial $Q$ of degree $d$ using explicit auxiliary functions, a technique introduced in number theory by C.J. Smyth to study the absolute trace of a totally positive algebraic integer \cite{Smyth1984}. To bound the roots of the polynomials for each of the degrees $d = 17, 18, 19$ and $20$, we use a version of the auxiliary function due to V. Flammang \cite{Flammang2016}, using which we establish that
\begin{table}[h]
\centering
\begin{tabular}{ccc}
\hline
\textbf{Degrees $d$} & \textbf{Trace $2d-3$} & \textbf{Intervals} \\ \hline
17                  & 31                    & (0,6.69)                         \\
18                  & 33                    & (0,7.11)                         \\
19                  & 35                    & (0,7.50)                         \\
20                  & 37                    & (0,7.86)                         \\\hline
\end{tabular}
\caption{Intervals containing the roots of $Q$ for degrees $d$}
\label{tab:intervalsRoots}
\end{table}

The second step is to propose bounds for the coefficients of the polynomials we are looking for. As explained in \cite{ElOtmaniMaulRhinSac-Epee2014} and \cite{ElOtmaniRhinSac-Epee2015}, this is achieved by uniformly drawing, in each interval, $d-1$ real numbers in $(0,4)$ and one real number greater than $4$. With these real numbers, we reconstruct the polynomial (with real coefficients) of which they are the roots..  This is repeated many times, enabling us to propose reliable bounds for the integer coefficients of the polynomials of interest, based on the location of the real coefficients of the polynomials obtained via our random draws.

The objective is then to identify monic integer polynomials $Q=\displaystyle\sum_{i=0}^d a_i x^i$ ($a_d = 1$), with all roots in the various intervals listed in Table~\ref{tab:intervalsRoots}, where the $a_i$'s are the unknown variables of the system we aim to solve. As mentioned above, polynomials $Q$ are required to have a single root greater than $4$, and all other roots in $(0,4)$. For an interval in Table~\ref{tab:intervalsRoots} corresponding to degree $d$, we first note that a polynomial $Q$ that vanishes $d$ times within this interval must alternate between strictly positive and strictly negative sign at points separating its roots. We leverage this observation as follows: we uniformly draw $d-2$ real numbers $\beta_i$ in interval $(0,4)$, and we add to these $d-2$ values the real numbers $0$, $4$, and the right-hand bound of the considered interval. We order these real numbers in ascending order, and we require the polynomial $Q$ to change sign at each of these $d+1$ real values, seen as roots separators. We must have $Q(0) < 0$ if $d$ is odd, $Q(0) >0$ if $d$ is even, and the sign of $Q$ must change when we move from one separator to the next.
This gives us a list of linear constraints that form the first part of our optimization problem. The system is completed by adding lower and upper bounds on the $a_i$'s. In \cite{ElOtmaniMaulRhinSac-Epee2014} and \cite{ElOtmaniRhinSac-Epee2015}, the libraries used required the minimization of an artificially introduced linear function, subject to a set of constraints, and the result of the minimization itself was of no interest, since discovering the $a_i$'s that satisfied the constraints was the sole true objective. With the newly experimented tools (the JuMP modeling language and the Gurobi optimization library), we were able to develop a much more powerful algorithm, focused on feasibility (that is, specific compliance with constraints) without the need to carry out the actual minimization of a linear function simply introduced for practical reasons.

\section{Implementation and results}

The improved algorithm has been fully implemented in the Julia programming language. Polynomial manipulations are performed using the \textit{Polynomials} package. For irreducibility tests, the \textit{Nemo} package is used.  The optimization problem is formulated within the framework of the \textit{JuMP} modeling language, while the search for a solution that respects the constraints imposed is carried out using the \textit{Gurobi} library via a wrapper maintained by the \textit{JuMP} community.

The algorithm launches a loop with the following three steps:
\begin{itemize}
  \item Randomly draw separators,
  \item Look for a vector of integers respecting the constraints,
  \item Test the irreducibility of the polynomial whose coefficients are the integers possibly found in the previous step.
\end{itemize}

Table~\ref{tab:totallyPositive17a19} contains two polynomials of degree $17$ and trace $31$, the second of which already appeared in \cite{ElOtmaniRhinSac-Epee2015}. It also includes four polynomials of degree $18$ and trace $33$, and fourteen polynomials of degree $19$ and trace  $35$. Table~\ref{tab:totallyPositive20} contains forty-five polynomials of degree $20$ and trace $37$.

Table~\ref{tab:SalemDegree34}, Table~\ref{tab:SalemDegree36}, Table~\ref{tab:SalemDegree38} and Table~\ref{tab:SalemDegree40} present the coefficients of the Salem polynomials of degree $2d$ obtained from polynomials of Table~\ref{tab:totallyPositive17a19} and Table~\ref{tab:totallyPositive20} by inverse change of variable. Only the coefficients necessary for the complete knowledge of each of the polynomials are shown in the tables, as we recall that these polynomials are reciprocal. In the last column of each table, the Salem numbers whose minimal polynomial appears on the same line are listed.

\section{Conclusion}
The refinement of our optimization algorithm has enabled us to achieve results that were out of reach with previous versions. We hope that the lists of Salem numbers with trace $-3$ are complete for degrees $34$, $36$, and possibly $38$, considering the number of times the algorithm has rediscovered each of them.

\begin{sidewaystable}
    \centering
    \begin{adjustbox}{width=\textheight} 
    \renewcommand{\arraystretch}{1.4}
    \scriptsize
    \begin{tabular}{|c|c|c|c|c|c|c|c|c|c|c|c|c|c|c|c|c|c|c|c|c|c|}
    \hline
    $x^{19}$ & $x^{18}$ & $x^{17}$ & $x^{16}$ & $x^{15}$ &$x^{14}$ & $x^{13}$ & $x^{12}$ & $x^{11}$ & $x^{10}$ & $x^9$ & $x^8$ & $x^7$ & $x^6$ & $x^5$ & $x^4$ & $x^3$ & $x^2$ & $x^1$ & $x^0$ \\
\hline
&  & 1 & -31 & 432 & -3583 & 19736 & -76285 & 213173 & -437297 & 662123 & -738327 & 600855 & -351222 & 143998 & -40022 & 7179 & -772 & 44 & -1 \\
&  & 1 & -31 & 433 & -3608 & 20013 & -78079 & 220717 & -458940 & 705459 & -799257 & 660596 & -391294 & 161786 & -44982 & 7979 & -837 & 46 & -1 \\
\hline
& 1 & -33 & 493 & -4418 & 26527 & -112896 & 351385 & -813619 & 1412249 & -1838442 & 1784606 & -1276751 & 661168 & -241606 & 60150 & -9713 & 946 & -49 & 1\\
& 1 & -33 & 493 & -4418 & 26526 & -112873 & 351153 & -812266 & 1407186 & -1825688 & 1762537 & -1250469 & 639916 & -230279 & 56361 & -8974 & 871 & -46 & 1\\
& 1 & -33 & 494 & -4445 & 26854 & -115243 & 362498 & -850194 & 1498144 & -1983921 & 1962162 & -1431230 & 754969 & -280049 & 70262 & -11292 & 1074 & -53 & 1\\
& 1 & -33 & 494 & -4445 & 26854 & -115243 & 362499 & -850212 & 1498282 & -1984511 & 1963710 & -1433817 & 757739 & -281916 & 71023 & -11468 & 1095 & -54 & 1\\
\hline
1 & -35 & 557 & -5344 & 34553 & -159439 & 542423 & -1386300 & 2687875 & -3966530 & 4445822 & -3758294 & 2367533 & -1091721 & 359274 & -81407 & 12063 & -1083 & 52 & -1\\
1 & -35 & 557 & -5344 & 34551 & -159390 & 541892 & -1382942 & 2674104 & -3928030 & 4370711 & -3655340 & 2268985 & -1027046 & 331107 & -73678 & 10832 & -984 & 49 & -1\\
1 & -35 & 557 & -5344 & 34551 & -159391 & 541914 & -1383154 & 2675283 & -3932228 & 4380750 & -3671796 & 2287527 & -1041244 & 338312 & -76006 & 11283 & -1031 & 51 & -1\\
1 & -35 & 558 & -5374 & 34956 & -162647 & 559300 & -1448266 & 2851355 & -4280967 & 4888760 & -4213778 & 2705933 & -1270169 & 424300 & -97142 & 14439 & -1283 & 59 & -1\\
1 & -35 & 558 & -5374 & 34955 & -162625 & 559086 & -1447051 & 2846876 & -4269675 & 4868764 & -4188594 & 2683303 & -1255726 & 417842 & -95176 & 14054 & -1240 & 57 & -1\\
1 & -35 & 558 & -5373 & 34930 & -162344 & 557208 & -1438734 & 2821126 & -4212363 & 4775880 & -4078746 & 2589161 & -1198202 & 393469 & -88329 & 12870 & -1130 & 53 & -1\\
1 & -35 & 558 & -5373 & 34931 & -162369 & 557485 & -1440530 & 2828706 & -4234284 & 4820423 & -4142921 & 2654512 & -1244549 & 415716 & -95214 & 14134 & -1248 & 57 & -1\\
1 & -35 & 558 & -5372 & 34904 & -162040 & 555095 & -1429011 & 2789911 & -4140583 & 4656415 & -3935052 & 2465780 & -1124387 & 363890 & -80862 & 11790 & -1053 & 51 & -1\\
1 & -35 & 558 & -5372 & 34904 & -162041 & 555117 & -1429223 & 2791089 & -4144766 & 4666361 & -3951199 & 2483725 & -1137903 & 370645 & -83028 & 12212 & -1098 & 53 & -1\\
1 & -35 & 558 & -5372 & 34904 & -162042 & 555138 & -1429415 & 2792093 & -4148084 & 4673591 & -3961718 & 2493842 & -1144106 & 372869 & -83387 & 12197 & -1086 & 52 & -1\\
1 & -35 & 558 & -5372 & 34905 & -162064 & 555350 & -1430594 & 2796295 & -4158182 & 4690409 & -3981464 & 2510414 & -1154167 & 377281 & -84729 & 12453 & -1112 & 53 & -1\\
1 & -35 & 559 & -5402 & 35310 & -165319 & 572714 & -1495491 & 2971138 & -4502372 & 5187028 & -4503873 & 2905886 & -1365213 & 454403 & -103264 & 15211 & -1340 & 61 & -1\\
1 & -35 & 560 & -5431 & 35691 & -168321 & 588542 & -1554487 & 3131239 & -4823342 & 5664341 & -5027914 & 3325004 & -1603764 & 547583 & -126895 & 18777 & -1616 & 69 & -1\\
1 & -35 & 560 & -5430 & 35664 & -167995 & 586218 & -1543608 & 3096048 & -4742710 & 5532312 & -4873839 & 3198484 & -1532318 & 520824 & -120639 & 17968 & -1574 & 69 & -1\\
    \hline
    \end{tabular}
    \end{adjustbox}
    \caption{Totally positive polynomials with degree $d = 17, 18, 19$ and trace $2d-3 = 31, 33, 35$, a single root $ > 4$, and all other roots in $(0,4)$.}
    \label{tab:totallyPositive17a19}
\end{sidewaystable}

\begin{sidewaystable}
    \centering
    \begin{adjustbox}{width=\textheight} 
    \renewcommand{\arraystretch}{1.4}
    \scriptsize
    \begin{tabular}{|c|c|c|c|c|c|c|c|c|c|c|c|c|c|c|c|c|c|c|c|c|c|c|}
    \hline
    $x^{20}$ & $x^{19}$ & $x^{18}$ & $x^{17}$ & $x^{16}$ & $x^{15}$ &$x^{14}$ & $x^{13}$ & $x^{12}$ & $x^{11}$ & $x^{10}$ & $x^9$ & $x^8$ & $x^7$ & $x^6$ & $x^5$ & $x^4$ & $x^3$ & $x^2$ & $x^1$ & $x^0$ \\
\hline
1 & -37 & 625 & -6399 & 44446 & -222060 & 825522 & -2329872 & 5048500 & -8439635 & 10881501 & -10768511 & 8104198 & -4574262 & 1899487 & -565429 & 116474 & -15800 & 1309 & -58 & 1\\
1 & -37 & 625 & -6398 & 44419 & -221730 & 823105 & -2318036 & 5007562 & -8336641 & 10690269 & -10505421 & 7837311 & -4377327 & 1796212 & -528252 & 107731 & -14551 & 1213 & -55 & 1\\
1 & -37 & 625 & -6398 & 44420 & -221758 & 823452 & -2320556 & 5019522 & -8375790 & 10781149 & -10656896 & 8018803 & -4532351 & 1888945 & -566002 & 117745 & -16168 & 1354 & -60 & 1\\
1 & -37 & 625 & -6398 & 44420 & -221756 & 823406 & -2320087 & 5016723 & -8364912 & 10752164 & -10602585 & 7946553 & -4464268 & 1844094 & -545867 & 111840 & -15112 & 1252 & -56 & 1\\
1 & -37 & 625 & -6398 & 44420 & -221754 & 823362 & -2319658 & 5014275 & -8355817 & 10729005 & -10561142 & 7893952 & -4417035 & 1814480 & -533222 & 108311 & -14511 & 1197 & -54 & 1\\
1 & -37 & 625 & -6397 & 44393 & -221427 & 821013 & -2308505 & 4977350 & -8268158 & 10577853 & -10371401 & 7721682 & -4305574 & 1764264 & -517945 & 105277 & -14129 & 1168 & -53 & 1\\
1 & -37 & 626 & -6429 & 44858 & -225492 & 844888 & -2408155 & 5282259 & -8963730 & 11769570 & -11904755 & 9193589 & -5345760 & 2294145 & -706514 & 150019 & -20708 & 1695 & -70 & 1\\
1 & -37 & 626 & -6428 & 44829 & -225112 & 841911 & -2392603 & 5225038 & -8811017 & 11469514 & -11468627 & 8726271 & -4980804 & 2090616 & -628026 & 130099 & -17632 & 1443 & -62 & 1\\
1 & -37 & 626 & -6428 & 44830 & -225137 & 842191 & -2394462 & 5233198 & -8836024 & 11524592 & -11557030 & 8829982 & -5069161 & 2144406 & -650788 & 136508 & -18750 & 1549 & -66 & 1\\
1 & -37 & 626 & -6428 & 44830 & -225136 & 842168 & -2394226 & 5231770 & -8830351 & 11509007 & -11526662 & 8787629 & -5027014 & 2114896 & -636621 & 132034 & -17881 & 1457 & -62 & 1\\
1 & -37 & 626 & -6428 & 44830 & -225136 & 842168 & -2394226 & 5231771 & -8830369 & 11509146 & -11527266 & 8789257 & -5029844 & 2118095 & -638940 & 133075 & -18153 & 1494 & -64 & 1\\
1 & -37 & 626 & -6428 & 44831 & -225164 & 842517 & -2396790 & 5244157 & -8871910 & 11608670 & -11699839 & 9006274 & -5226294 & 2243942 & -694391 & 149148 & -21004 & 1765 & -74 & 1\\
1 & -37 & 626 & -6428 & 44831 & -225163 & 842493 & -2396532 & 5242515 & -8865020 & 11588579 & -11658031 & 8943510 & -5158379 & 2191598 & -666355 & 139149 & -18798 & 1503 & -62 & 1\\
1 & -37 & 626 & -6428 & 44831 & -225162 & 842470 & -2396297 & 5241107 & -8859521 & 11573861 & -11630406 & 8906932 & -5124480 & 2170088 & -657381 & 136857 & -18486 & 1487 & -62 & 1\\
1 & -37 & 626 & -6427 & 44802 & -224783 & 839517 & -2381002 & 5185507 & -8713504 & 11292862 & -11232573 & 8494427 & -4815249 & 2006263 & -598222 & 123099 & -16607 & 1359 & -59 & 1\\
1 & -37 & 626 & -6427 & 44802 & -224781 & 839470 & -2380512 & 5182516 & -8701623 & 11260573 & -11171114 & 8411920 & -4737543 & 1955761 & -576212 & 116949 & -15583 & 1270 & -56 & 1\\
1 & -37 & 626 & -6427 & 44803 & -224810 & 839841 & -2383289 & 5196099 & -8747502 & 11370755 & -11361866 & 8650639 & -4952078 & 2092024 & -635760 & 134120 & -18639 & 1568 & -68 & 1\\
1 & -37 & 626 & -6427 & 44803 & -224809 & 839817 & -2383031 & 5194458 & -8740630 & 11350803 & -11320664 & 8589526 & -4887098 & 2043125 & -610360 & 125384 & -16782 & 1353 & -58 & 1\\
1 & -37 & 626 & -6427 & 44803 & -224809 & 839819 & -2383072 & 5194825 & -8742517 & 11356972 & -11334072 & 8609206 & -4906517 & 2055705 & -615459 & 126566 & -16912 & 1357 & -58 & 1\\
1 & -37 & 626 & -6427 & 44803 & -224807 & 839772 & -2382584 & 5191872 & -8730949 & 11326157 & -11277021 & 8535488 & -4840735 & 2016215 & -600243 & 123098 & -16512 & 1341 & -58 & 1\\
1 & -37 & 627 & -6457 & 45212 & -228165 & 858322 & -2455549 & 5402809 & -9187016 & 12069406 & -12190247 & 9375898 & -5412124 & 2297027 & -696834 & 145282 & -19678 & 1589 & -66 & 1\\
1 & -37 & 627 & -6457 & 45213 & -228192 & 858650 & -2457919 & 5414155 & -9224963 & 12160520 & -12349205 & 9577625 & -5596843 & 2416972 & -750527 & 161153 & -22574 & 1879 & -78 & 1\\
1 & -37 & 627 & -6457 & 45213 & -228191 & 858625 & -2457641 & 5412339 & -9217208 & 12137715 & -12301797 & 9507190 & -5522069 & 2360827 & -721344 & 151042 & -20390 & 1621 & -66 & 1\\
1 & -37 & 627 & -6457 & 45213 & -228191 & 858625 & -2457641 & 5412340 & -9217226 & 12137855 & -12302416 & 9508912 & -5525218 & 2364663 & -724426 & 152619 & -20868 & 1695 & -70 & 1\\
1 & -37 & 627 & -6456 & 45184 & -227810 & 855625 & -2441857 & 5353767 & -9059466 & 11825157 & -11844671 & 9016474 & -5140987 & 2152040 & -643778 & 132678 & -17884 & 1455 & -62 & 1\\
1 & -37 & 627 & -6456 & 45184 & -227810 & 855625 & -2441856 & 5353748 & -9059310 & 11824430 & -11842546 & 9012408 & -5135825 & 2147726 & -641477 & 131945 & -17762 & 1447 & -62 & 1\\
1 & -37 & 627 & -6456 & 45184 & -227810 & 855625 & -2441856 & 5353748 & -9059309 & 11824413 & -11842425 & 9011936 & -5134715 & 2146096 & -639982 & 131115 & -17503 & 1408 & -60 & 1\\
1 & -37 & 627 & -6456 & 45184 & -227810 & 855625 & -2441856 & 5353748 & -9059309 & 11824414 & -11842438 & 9012005 & -5134910 & 2146416 & -640293 & 131289 & -17554 & 1414 & -60 & 1\\
1 & -37 & 627 & -6456 & 45184 & -227810 & 855625 & -2441856 & 5353748 & -9059309 & 11824414 & -11842438 & 9012005 & -5134909 & 2146409 & -640275 & 131268 & -17543 & 1412 & -60 & 1\\
1 & -37 & 627 & -6456 & 45184 & -227810 & 855626 & -2441878 & 5353960 & -9060487 & 11828596 & -11852371 & 9028084 & -5152669 & 2159643 & -646788 & 133322 & -17940 & 1455 & -62 & 1\\
1 & -37 & 627 & -6456 & 45185 & -227838 & 855974 & -2444417 & 5366077 & -9100380 & 11921723 & -12008396 & 9215903 & -5313756 & 2256284 & -686103 & 143642 & -19552 & 1585 & -66 & 1\\
1 & -37 & 627 & -6456 & 45185 & -227838 & 855975 & -2444438 & 5366272 & -9101438 & 11925456 & -12017410 & 9231150 & -5331914 & 2271354 & -694588 & 146725 & -20212 & 1654 & -68 & 1\\
1 & -37 & 627 & -6456 & 45185 & -227837 & 855951 & -2444183 & 5364688 & -9095036 & 11907718 & -11982813 & 9183119 & -5284481 & 2238431 & -678944 & 141878 & -19318 & 1573 & -66 & 1\\
1 & -37 & 627 & -6456 & 45185 & -227837 & 855951 & -2444181 & 5364649 & -9094706 & 11906123 & -11977932 & 9173209 & -5270879 & 2225814 & -671175 & 138825 & -18606 & 1487 & -62 & 1\\
1 & -37 & 627 & -6456 & 45185 & -227837 & 855952 & -2444206 & 5364922 & -9096425 & 11913062 & -11996815 & 9208666 & -5317092 & 2267276 & -696202 & 148578 & -20895 & 1772 & -76 & 1\\
1 & -37 & 627 & -6456 & 45185 & -227837 & 855952 & -2444204 & 5364880 & -9096041 & 11911051 & -11990136 & 9193950 & -5295224 & 2245421 & -681799 & 142556 & -19400 & 1577 & -66 & 1\\
1 & -37 & 627 & -6456 & 45185 & -227835 & 855904 & -2443694 & 5361716 & -9083310 & 11876144 & -11923421 & 9104529 & -5211706 & 2192000 & -659041 & 136344 & -18369 & 1481 & -62 & 1\\
1 & -37 & 627 & -6456 & 45185 & -227835 & 855904 & -2443692 & 5361676 & -9082964 & 11874441 & -11918141 & 9093749 & -5196994 & 2178655 & -651191 & 133485 & -17773 & 1420 & -60 & 1\\
1 & -37 & 628 & -6486 & 45592 & -231142 & 873873 & -2512751 & 5555367 & -9486362 & 12503537 & -12654056 & 9737224 & -5614086 & 2376431 & -718438 & 149359 & -20224 & 1638 & -68 & 1\\
1 & -37 & 628 & -6486 & 45592 & -231142 & 873873 & -2512750 & 5555346 & -9486170 & 12502532 & -12650724 & 9729912 & -5603303 & 2365794 & -711572 & 146576 & -19562 & 1556 & -64 & 1\\
1 & -37 & 628 & -6486 & 45592 & -231142 & 873874 & -2512772 & 5555560 & -9487383 & 12506975 & -12661742 & 9748752 & -5625507 & 2383541 & -720881 & 149612 & -20128 & 1609 & -66 & 1\\
1 & -37 & 628 & -6486 & 45593 & -231169 & 874201 & -2515119 & 5566671 & -9523924 & 12592623 & -12806213 & 9923750 & -5775760 & 2472712 & -755964 & 158148 & -21258 & 1669 & -66 & 1\\
1 & -37 & 628 & -6486 & 45593 & -231168 & 874176 & -2514842 & 5564877 & -9516380 & 12570979 & -12762865 & 9862764 & -5715893 & 2432511 & -738161 & 153272 & -20528 & 1625 & -66 & 1\\
1 & -37 & 628 & -6486 & 45593 & -231167 & 874152 & -2514587 & 5563295 & -9510015 & 12553535 & -12729583 & 9818381 & -5674984 & 2407163 & -728125 & 150950 & -20263 & 1615 & -66 & 1\\
1 & -37 & 628 & -6485 & 45564 & -230788 & 871199 & -2499292 & 5507696 & -9364014 & 12272639 & -12332084 & 9406401 & -5365891 & 2242486 & -667539 & 136173 & -18026 & 1430 & -60 & 1\\
\hline
    \end{tabular}
    \end{adjustbox}
    \caption{Totally positive polynomials with degree $d = 20$ and trace $2d-3 = 37$, a single root $ > 4$, and all other roots in $(0,4)$.}
    \label{tab:totallyPositive20}
\end{sidewaystable}

\vspace{5pt}
\begin{table}[ht]
\centering
\begin{adjustbox}{width=\textwidth}
\renewcommand{\arraystretch}{1.4}
\begin{tabular}{|c|c|c|c|c|c|c|c|c|c|c|c|c|c|c|c|c|c|c|}
\hline
$x^{34}$ & $x^{33}$ & $x^{32}$ & $x^{31}$ & $x^{30}$ & $x^{29}$ &$x^{28}$ & $x^{27}$ & $x^{26}$ & $x^{25}$ & $x^{24}$ & $x^{23}$ & $x^{22}$ & $x^{21}$ & $x^{20}$ & $x^{19}$ & $x^{18}$ & $x^{17}$ & Salem numbers\\
\hline
1& 3 &1 &-15& -52 &-107 &-167 &-220 &-266 &-317 &-383 &-462 &-540 &-603 &-646 &-673 &-688 &-693 & 2.9718375039\\
1 &3 &2& -10 &-40 &-89 &-149 &-208 &-257& -293& -315& -322& -311 &-281 &-237& -191 &-156& -143 & 2.7616448085 \\
\hline
\end{tabular}
\end{adjustbox}
\caption{Salem numbers of degree $34$ and trace $-3$.} 
\label{tab:SalemDegree34}
\end{table}

\vspace{5pt}
\begin{table}[ht]
\centering
\begin{adjustbox}{width=\textwidth}
\renewcommand{\arraystretch}{1.4}
\begin{tabular}{|c|c|c|c|c|c|c|c|c|c|c|c|c|c|c|c|c|c|c|c|}
\hline
$x^{36}$ & $x^{35}$ & $x^{34}$ & $x^{33}$ & $x^{32}$ & $x^{31}$ & $x^{30}$ & $x^{29}$ &$x^{28}$ & $x^{27}$ & $x^{26}$ & $x^{25}$ & $x^{24}$ & $x^{23}$ & $x^{22}$ & $x^{21}$ & $x^{20}$ & $x^{19}$ & $x^{18}$ & Salem numbers\\
\hline
1 &3 &1 &-15& -52& -106& -161& -201& -224& -245& -282& -340& -404& -452& -473& -474& -469& -467& -467& 2.9585975515\\
1 &3 &1 &-15& -53& -111& -173& -217& -228& -206& -161 &-103 &-34& 47& 136 &221 &289 &331 &345 & 2.9778824718\\
1 &3 &2& -10 &-39& -82 &-124& -146 &-136 &-96& -40& 12& 45& 55& 50& 43& 43& 48& 51& 2.6852347917\\
1& 3 &2& -10 &-39& -82 &-123& -140 &-118& -60& 14 &77& 111& 113& 92& 61 &32& 12 &5 & 2.6779781031\\
\hline
\end{tabular}
\end{adjustbox}
\caption{Salem numbers of degree $36$ and trace $-3$.} 
\label{tab:SalemDegree36}
\end{table}

\begin{table}[ht]
\centering
\begin{adjustbox}{width=\textwidth}
\renewcommand{\arraystretch}{1.4}
\begin{tabular}{|c|c|c|c|c|c|c|c|c|c|c|c|c|c|c|c|c|c|c|c|c|}
\hline
$x^{38}$ & $x^{37}$ & $x^{36}$ & $x^{35}$ & $x^{34}$ & $x^{33}$ & $x^{32}$ & $x^{31}$ & $x^{30}$ & $x^{29}$ &$x^{28}$ & $x^{27}$ & $x^{26}$ & $x^{25}$ & $x^{24}$ & $x^{23}$ & $x^{22}$ & $x^{21}$ & $x^{20}$ & $x^{19}$ &  Salem numbers\\
\hline
1& 3& 0& -20 &-63& -118& -161& -176& -173& -183& -228& -298& -357& -376& -360& -345& -364& -419& -479 &-505& 3.1314069583\\
1& 3 &0 &-20& -65& -129& -190& -222& -213& -174& -128& -94 &-77& -70& -61& -42 &-11& 26& 57& 69& 3.1753281563\\
1& 3& 0& -20& -65& -130& -196& -240& -248& -221& -168& -100& -27& 43& 107& 166& 223& 276& 316& 331& 3.1837638764\\
1& 3& 1& -16& -59& -132& -227& -328& -421& -502& -575& -648& -726& -807& -882& -940& -972& -979& -973& -969& 3.1073115409\\
1& 3& 1& -16& -60& -140& -260& -421& -623& -864& -1135& -1422& -1710& -1986& -2238& -2455& -2626& -2745& -2813& -2835& 3.1627585984\\
1& 3& 1& -15& -53& -113& -185& -254& -308& -343& -360& -364& -363& -367& -384& -417& -460& -502& -532& -543& 3.0026153431\\
1& 3& 1& -15& -52& -108& -173& -238& -302& -369& -436& -491& -522& -528& -518& -502& -482& -458& -436& -427& 2.9830177489\\
1& 3& 1& -14& -47& -93& -136& -157& -146& -109& -62& -22& 2& 12& 18& 29& 49& 74& 95& 103& 2.8626146783\\
1& 3& 1& -14& -47& -94& -142& -175& -182& -163& -126& -83& -45& -18& 0& 15& 33& 54& 72& 79& 2.8769732691\\
1& 3& 1& -14& -47& -95& -149& -199& -237& -260& -270& -274& -280& -294& -318& -350& -382& -406& -419& -423& 2.8963426610\\
1& 3& 1& -14& -46& -87& -118& -122& -98& -62& -31& -10& 9& 35& 67& 93& 103& 97& 85& 79& 2.8205243729\\
1& 3& 2& -10& -40& -88& -143& -189& -214& -217& -208& -203& -214& -244& -285& -323& -346& -352& -348& -345& 2.7432751935\\
1& 3& 3& -5& -26& -60& -100& -137& -164& -177& -173& -150& -107& -47& 25& 104& 185& 259& 313& 333& 2.4399913855\\
1& 3& 3& -4& -21& -48& -81& -114& -142& -163& -176& -183& -188& -195& -206& -221& -236& -247& -253& -255& 2.3319842011\\
\hline
\end{tabular}
\end{adjustbox}
\caption{Salem numbers of degree $38$ and trace $-3$.} 
\label{tab:SalemDegree38}
\end{table}

\begin{table}[ht]
\centering
\begin{adjustbox}{width=\textwidth}
\renewcommand{\arraystretch}{1.4}
\begin{tabular}{|c|c|c|c|c|c|c|c|c|c|c|c|c|c|c|c|c|c|c|c|c|c|}
\hline
$x^{40}$ & $x^{39}$ & $x^{38}$ & $x^{37}$ & $x^{36}$ & $x^{35}$ & $x^{34}$ & $x^{33}$ & $x^{32}$ & $x^{31}$ & $x^{30}$ & $x^{29}$ &$x^{28}$ & $x^{27}$ & $x^{26}$ & $x^{25}$ & $x^{24}$ & $x^{23}$ & $x^{22}$ & $x^{21}$ & $x^{20}$ &  Salem numbers\\
\hline
1 & 3 & -1 & -30 & -112 & -274 & -527 & -857 & -1227 & -1587 & -1883 & -2068 & -2110 & -1998 & -1744 & -1382 & -960 & -535 & -168 & 82 & 171 & 3.8007708538 \\
1 & 3 & -1 & -29 & -105 & -247 & -452 & -691 & -919 & -1095 & -1198 & -1234 & -1229 & -1217 & -1226 & -1268 & -1335 & -1408 & -1468 & -1505 & -1517 & 3.7206401214 \\
1 & 3 & -1 & -29 & -104 & -243 & -449 & -715 & -1033 & -1396 & -1787 & -2172 & -2507 & -2759 & -2920 & -3003 & -3025 & -3001 & -2949 & -2897 & -2875 & 3.7193813271 \\
1 & 3 & -1 & -29 & -104 & -241 & -435 & -664 & -906 & -1155 & -1418 & -1699 & -1983 & -2239 & -2440 & -2584 & -2693 & -2794 & -2894 & -2974 & -3005 & 3.7085516704 \\
1 & 3 & -1 & -29 & -104 & -239 & -419 & -597 & -712 & -721 & -623 & -463 & -308 & -210 & -177 & -172 & -138 & -40 & 107 & 242 & 297 & 3.6933302502 \\
1 & 3 & -1 & -28 & -97 & -215 & -366 & -515 & -625 & -674 & -658 & -586 & -474 & -346 & -232 & -157 & -125 & -120 & -121 & -119 & -117 & 3.6222259269 \\
1 & 3 & 0 & -24 & -90 & -212 & -384 & -580 & -766 & -913 & -998 & -999 & -891 & -655 & -293 & 163 & 658 & 1127 & 1510 & 1760 & 1847 & 3.5263915198 \\
1 & 3 & 0 & -23 & -85 & -199 & -361 & -548 & -725 & -856 & -911 & -872 & -735 & -512 & -229 & 80 & 383 & 653 & 868 & 1008 & 1057 & 3.4775939743 \\
1 & 3 & 0 & -23 & -84 & -192 & -335 & -482 & -599 & -667 & -685 & -661 & -598 & -492 & -343 & -168 & 2 & 136 & 220 & 260 & 271 & 3.4505894081 \\
1 & 3 & 0 & -23 & -84 & -191 & -328 & -455 & -525 & -509 & -411 & -267 & -125 & -21 & 37 & 68 & 104 & 166 & 249 & 323 & 353 & 3.4416559110 \\
1 & 3 & 0 & -23 & -84 & -191 & -328 & -455 & -524 & -503 & -392 & -225 & -53 & 79 & 151 & 173 & 176 & 189 & 221 & 257 & 273 & 3.4415063642 \\
1 & 3 & 0 & -23 & -83 & -187 & -323 & -467 & -600 & -720 & -835 & -949 & -1054 & -1139 & -1203 & -1257 & -1311 & -1365 & -1410 & -1438 & -1447 & 3.4387135979 \\
1 & 3 & 0 & -23 & -83 & -186 & -317 & -446 & -546 & -607 & -632 & -624 & -576 & -477 & -326 & -141 & 48 & 211 & 330 & 400 & 423 & 3.4313736405 \\
1 & 3 & 0 & -23 & -83 & -185 & -310 & -420 & -478 & -468 & -399 & -297 & -192 & -113 & -87 & -135 & -260 & -440 & -629 & -772 & -825 & 3.4226176184 \\
1 & 3 & 0 & -22 & -78 & -173 & -293 & -408 & -486 & -509 & -477 & -402 & -297 & -173 & -42 & 80 & 178 & 243 & 277 & 290 & 293 & 3.3753589936 \\
1 & 3 & 0 & -22 & -78 & -171 & -280 & -364 & -383 & -324 & -213 & -106 & -58 & -91 & -178 & -255 & -255 & -151 & 23 & 185 & 251 & 3.3596154717 \\
1 & 3 & 0 & -22 & -77 & -168 & -283 & -408 & -544 & -709 & -918 & -1164 & -1417 & -1646 & -1840 & -2009 & -2165 & -2308 & -2425 & -2501 & -2527 & 3.3686279815 \\
1 & 3 & 0 & -22 & -77 & -167 & -277 & -387 & -489 & -590 & -696 & -799 & -880 & -927 & -945 & -948 & -942 & -924 & -893 & -861 & -847 & 3.3601726633 \\
1 & 3 & 0 & -22 & -77 & -167 & -275 & -372 & -432 & -442 & -399 & -306 & -172 & -17 & 130 & 242 & 309 & 339 & 349 & 352 & 353 & 3.3538404471 \\
1 & 3 & 0 & -22 & -77 & -165 & -262 & -330 & -343 & -308 & -259 & -229 & -220 & -202 & -142 & -38 & 76 & 155 & 180 & 170 & 161 & 3.3390688183 \\
1 & 3 & 1 & -16 & -58 & -125 & -203 & -273 & -324 & -358 & -382 & -400 & -410 & -410 & -403 & -396 & -393 & -394 & -396 & -397 & -397 & 3.0674618430 \\
1 & 3 & 1 & -16 & -57 & -120 & -189 & -244 & -274 & -282 & -276 & -260 & -231 & -187 & -134 & -84 & -44 & -13 & 13 & 33 & 41 & 3.0404935832 \\
1 & 3 & 1 & -16 & -57 & -119 & -184 & -231 & -252 & -256 & -256 & -256 & -248 & -223 & -181 & -131 & -81 & -35 & 5 & 34 & 45 & 3.0329139322 \\
1 & 3 & 1 & -16 & -57 & -119 & -184 & -231 & -251 & -250 & -236 & -209 & -162 & -93 & -11 & 70 & 143 & 208 & 265 & 307 & 323 & 3.0324565729 \\
1 & 3 & 1 & -15 & -52 & -105 & -154 & -176 & -164 & -137 & -128 & -161 & -232 & -312 & -368 & -386 & -375 & -356 & -344 & -341 & -341 & 2.9406722181 \\
1 & 3 & 1 & -15 & -52 & -105 & -154 & -175 & -157 & -112 & -67 & -46 & -52 & -67 & -69 & -52 & -29 & -20 & -32 & -53 & -63 & 2.9384433127 \\
1 & 3 & 1 & -15 & -52 & -105 & -154 & -175 & -157 & -111 & -62 & -34 & -36 & -61 & -94 & -122 & -137 & -138 & -129 & -118 & -113 & 2.9383434641 \\
1 & 3 & 1 & -15 & -52 & -105 & -154 & -175 & -157 & -111 & -61 & -27 & -11 & -1 & 13 & 26 & 23 & -7 & -57 & -104 & -123 & 2.9382595870 \\
1 & 3 & 1 & -15 & -52 & -105 & -154 & -175 & -157 & -111 & -61 & -27 & -11 & 0 & 20 & 51 & 84 & 107 & 116 & 116 & 115 & 2.9382561494 \\
1 & 3 & 1 & -15 & -52 & -105 & -153 & -169 & -139 & -75 & -8 & 30 & 26 & -7 & -39 & -42 & -2 & 75 & 166 & 239 & 267 & 2.9339655654 \\
1 & 3 & 1 & -15 & -51 & -101 & -149 & -184 & -212 & -251 & -311 & -384 & -453 & -512 & -571 & -640 & -709 & -754 & -761 & -744 & -733 & 2.9358570973 \\
1 & 3 & 1 & -15 & -51 & -101 & -148 & -177 & -185 & -178 & -159 & -126 & -78 & -20 & 41 & 100 & 153 & 192 & 211 & 214 & 213 & 2.9281342716 \\
1 & 3 & 1 & -15 & -51 & -100 & -142 & -159 & -151 & -136 & -131 & -138 & -147 & -153 & -164 & -191 & -230 & -263 & -275 & -270 & -265 & 2.9161201862 \\
1 & 3 & 1 & -15 & -51 & -100 & -142 & -157 & -138 & -92 & -30 & 37 & 97 & 136 & 144 & 120 & 72 & 10 & -54 & -104 & -123 & 2.9120722906 \\
1 & 3 & 1 & -15 & -51 & -100 & -141 & -154 & -137 & -108 & -86 & -73 & -52 & -7 & 58 & 120 & 158 & 168 & 162 & 154 & 151 & 2.9121557127 \\
1 & 3 & 1 & -15 & -51 & -100 & -141 & -152 & -127 & -82 & -41 & -18 & -8 & 4 & 28 & 61 & 94 & 120 & 138 & 149 & 153 & 2.9098558542 \\
1 & 3 & 1 & -15 & -51 & -98 & -129 & -116 & -55 & 25 & 82 & 92 & 67 & 40 & 34 & 44 & 48 & 29 & -9 & -46 & -61 & 2.8817787817 \\
1 & 3 & 1 & -15 & -51 & -98 & -129 & -114 & -43 & 61 & 151 & 178 & 117 & -18 & -174 & -278 & -269 & -135 & 73 & 260 & 335 & 2.8785434334 \\
1 & 3 & 2 & -9 & -34 & -71 & -113 & -154 & -192 & -227 & -256 & -274 & -276 & -260 & -227 & -181 & -128 & -77 & -36 & -10 & -1 & 2.6316862080 \\
1 & 3 & 2 & -9 & -34 & -71 & -113 & -153 & -187 & -214 & -233 & -244 & -248 & -247 & -243 & -239 & -238 & -243 & -253 & -263 & -267 & 2.6296111462 \\
1 & 3 & 2 & -9 & -34 & -71 & -112 & -147 & -167 & -165 & -137 & -89 & -38 & -4 & 2 & -19 & -58 & -106 & -154 & -191 & -205 & 2.6188371017 \\
1 & 3 & 2 & -9 & -33 & -66 & -99 & -123 & -132 & -124 & -100 & -66 & -32 & -8 & 2 & 0 & -9 & -22 & -36 & -47 & -51 & 2.5752377703 \\
1 & 3 & 2 & -9 & -33 & -65 & -94 & -111 & -114 & -106 & -89 & -65 & -38 & -14 & 4 & 18 & 33 & 50 & 66 & 77 & 81 & 2.5602938026 \\
1 & 3 & 2 & -9 & -33 & -64 & -88 & -93 & -78 & -53 & -29 & -11 & 3 & 15 & 22 & 18 & 3 & -16 & -30 & -36 & -37 & 2.5347241894 \\
1 & 3 & 2 & -8 & -28 & -52 & -71 & -81 & -85 & -86 & -80 & -62 & -33 & -1 & 26 & 46 & 63 & 78 & 89 & 95 & 97 & 2.4352709937 \\
\hline
\end{tabular}
\end{adjustbox}
\caption{Salem numbers of degree $40$ and trace $-3$.} 
\label{tab:SalemDegree40}
\end{table}

\end{document}